\documentclass{elsarticle}
\usepackage[margin=1.25 in]{geometry}

\usepackage{amssymb}
\usepackage{hyperref}
\usepackage{amsfonts}
\usepackage{epsfig, amsmath}
\usepackage{tikz,ifthen,calc}
\usepackage{graphicx}
\usepackage{float}
\usepackage{caption}
\usepackage[labelfont=rm]{subcaption}
\usepackage{rotating}
\usetikzlibrary{patterns,snakes}
\newdefinition{rmk}{Remark}[section]
\newproof{pf}{Proof}
\newtheorem{thm}{Theorem}[section]

\title{The tilings of deficient squares by ribbon $L$-tetrominoes are diagonally cracked}

\author{Viorel Nitica}
\ead{vnitica@wcupa.edu}
\address[vn]{Department of Mathematics, West Chester University, West Chester, PA 19380}

\begin{document}
\begin{abstract} We consider tilings of deficient rectangles by the set $\mathcal{T}_4$ of ribbon $L$-tetrominoes. A tiling exists if and only if the rectangle is a square of odd side. The missing cell has to be on the main NW--SE diagonal, in an odd position if the square is $(4m+1)\times (4m+1)$ and in an even position if the square is $(4m+3)\times (4m+3)$. The majority of the tiles in a tiling follow the rectangular pattern, that is, are paired and each pair tiles a $2\times 4$ rectangle. The tiles in an irregular position together with the missing cell form a NW--SE diagonal crack. The crack is located in a thin region symmetric about the diagonal, made out of a sequence of $3\times 3$ squares that overlap over one of the corner cells. The crack divides the square in two parts of equal area. The number of tilings of a $(4m+1)\times (4m+1)$ deficient square by  $\mathcal{T}_4$ is equal to the number of tilings by dominoes of a $2m\times 2m$ square. The number of tilings of a $(4m+3)\times (4m+3)$ deficient square by  $\mathcal{T}_4$ is twice the number of tilings by dominoes of a $(2m+1)\times (2m+1)$ deficient square, with the missing cell placed on the main diagonal. In both cases the counting is realized by an explicit function which is a bijection in the first case and a double cover in the second. The crack in a square naturally propagates to a crack in a larger square. 

If an extra $2\times 2$ tile is added to $\mathcal{T}_4$, we call the new tile set $\mathcal{T}_4^+$. A tiling of a deficient rectangle by $\mathcal{T}_4^+$ exists if and only if the rectangle is a square of odd side. The missing cell has to be on the main NW--SE diagonal, in an odd position if the square is $(4m+1)\times (4m+1)$ and in an even position if the square is $(4m+3)\times (4m+3)$. The majority of the tiles in a tiling follow the rectangular pattern, that is, are either paired tetrominoes and each pair tiles a $2\times 4$ rectangle, or are $2\times 2$ squares. The tiles in an irregular position together with the missing cell form a NW--SE diagonal crack. The crack is located in a thin region symmetric about the diagonal, made out of a sequence of $3\times 3$ squares that overlap over one of the corner cells. The number of tilings of a $(4m+1)\times (4m+1)$ deficient square by  $\mathcal{T}_4^+$ is greater than the number of tilings by dominoes and monomers of a $2m\times 2m$ square. The number of tilings of a $(4m+3)\times (4m+3)$ deficient square by  $\mathcal{T}_4^+$ is greater than twice the number of tilings by dominoes and monomers of a $(2m+1)\times (2m+1)$ deficient square, with the missing cell placed on the main diagonal.The crack in a square naturally propagates to a crack in a larger square.

The results can be easily extended and allow for modeling of cracks in more irregular, even with extra holes, regions. 
\end{abstract}
	
\maketitle

\section{Introduction}\label{intro} We study tilings of deficient rectangles placed in a square lattice by tile sets consisting of polyominoes. A $1\times 1$ square in the lattice is called a {\em cell} or {\em monomer}. We call a rectangle {\em deficient} if a cell is missing. Our tile sets consist either only of ribbon $L$-tetrominoes ($\mathcal{T}_4$) or of ribbon $L$-tetrominoes and a $2\times 2$ square ($\mathcal{T}_4^+$). A \emph{tiling} of a region is a covering without overlapping of tiles. Only translations of the tiles are allowed in a tiling. We will call {\em $2$-square} a $2\times 2$ square that has the coordinates of all vertices even.

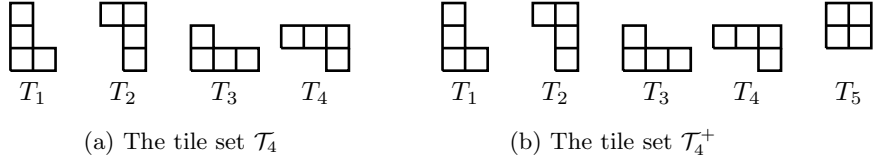
\begin{figure}[h]
\center
\begin{subfigure}{0.3\textwidth}
\begin{tikzpicture}[scale=.3]

\draw [line width = 1](0, 0) -- (2, 0) -- (2, 1) -- (1, 1) -- (1, 3) -- (0, 3) -- (0, 0);
\draw [line width = 1](0, 1)-- (1,1);
\draw [line width = 1](0, 2)-- (1,2);
\draw [line width = 1](1, 0)--(1,1);
\node at (1,-1){$T_1$};

\draw [line width = 1](5,1)--(6,1);
\draw [line width = 1](5,2)--(6,2);
\draw [line width = 1](5,3)--(5,2);
\draw [line width =1](5,0)--(6,0)--(6,3)--(4,3)--(4,2)--(5,2)--(5,0);
\node at (5, -1){$T_2$};

 \draw [line width = 1] (8,0)--(11,0)--(11,1)--(9,1)--(9,2)--(8,2)--(8,0);
 \draw [line width = 1] (8,1)--(9,1);
 \draw [line width = 1] (9,0)--(9,1);
 \draw [line width = 1] (10,0)--(10,1);

 \draw [line width = 1] (13,1)--(13,2);
 \draw [line width = 1] (15,1)--(15,0)--(14, 0)--(14,2);
 \draw [line width = 1] (12,1)--(15,1)--(15,2)--(12,2)--(12,1);
 \node at (9.5,-1){$T_3$};
 \node at (13.5, -1) {$T_4$};
\end{tikzpicture}
\caption{The tile set $\mathcal{T}_4$}
\end{subfigure}
~~~~~~~~
\begin{subfigure}{0.3\textwidth}
\begin{tikzpicture}[scale=.3]

\draw [line width = 1](0, 0) -- (2, 0) -- (2, 1) -- (1, 1) -- (1, 3) -- (0, 3) -- (0, 0);
\draw [line width = 1](0, 1)-- (1,1);
\draw [line width = 1](0, 2)-- (1,2);
\draw [line width = 1](1, 0)--(1,1);
\node at (1,-1){$T_1$};

\draw [line width = 1](5,1)--(6,1);
\draw [line width = 1](5,2)--(6,2);
\draw [line width = 1](5,3)--(5,2);
\draw [line width =1](5,0)--(6,0)--(6,3)--(4,3)--(4,2)--(5,2)--(5,0);
\node at (5, -1){$T_2$};

 \draw [line width = 1] (8,0)--(11,0)--(11,1)--(9,1)--(9,2)--(8,2)--(8,0);
 \draw [line width = 1] (8,1)--(9,1);
 \draw [line width = 1] (9,0)--(9,1);
 \draw [line width = 1] (10,0)--(10,1);

 \draw [line width = 1] (13,1)--(13,2);
 \draw [line width = 1] (15,1)--(15,0)--(14, 0)--(14,2);
 \draw [line width = 1] (12,1)--(15,1)--(15,2)--(12,2)--(12,1);
 \node at (9.5,-1){$T_3$};
 \node at (13.5, -1) {$T_4$};

\draw [line width = 1] (17,1)--(19,1)--(19,3)--(17,3)--(17,1);
 \draw [line width = 1] (18,1)--(18,3);
 \draw [line width = 1] (17,2)--(19,2);
 \node at (18,-1){$T_5$};

\end{tikzpicture}
\caption{The tile set $\mathcal{T}_4^+$}
\end{subfigure}
\caption{Tiles sets.}
\end{figure}

%%%%%%%%%%%%%%%%%%%%%%%%%%%%%%%%%%%%%%%%%%%%%%%%%%%%%%%%%%%%%%%%%%

\begin{figure}[h]
\center
\begin{subfigure}{0.3\textwidth}
\begin{tikzpicture}[scale=.3]
\draw (0,0) grid (11,11);
\draw [fill = darkgray] (3,7)--(4,7)--(4,8)--(3,8)--(3,7);

\draw [line width = 1](0, 0) -- (4, 0) -- (4,4) -- (0,4)--(0, 0);
\draw [line width = 1](0,2)-- (4,2);

\draw [line width = 1](4, 0) -- (8, 0) -- (8,4) -- (4,4)--(4, 0);
\draw [line width = 1](4,2)-- (8,2);

\draw [line width = 1] (2,4)--(6,4)--(6,6)--(2,6)--(2,4);
\draw [line width = 1] (9,1)--(11,1)--(11,5)--(9,5)--(9,1);
\draw [line width = 1] (9,7)--(11,7)--(11,11)--(9,11)--(9,7);
\draw [line width = 1] (7,7)--(9,7)--(9,11)--(7,11)--(7,7);
\draw [line width = 1] (5,7)--(7,7)--(7,11)--(5,11)--(5,7);
\draw [line width = 1] (11,5)--(11,7);
\draw [line width = 1] (1,11)--(5,11);

\draw [line width = 1] (0,4)--(2,4)--(2,8)--(0,8)--(0,4);

\draw [line width = 1, fill=lightgray] (0,8)--(2,8)--(2,9)--(1,9)--(1,11)--(0,11)--(0,8);
\draw [line width = 1, fill=lightgray] (2,6)--(4,6)--(4,7)--(3,7)--(3,9)--(2,9)--(2,6);
\draw [line width = 1, fill=lightgray] (4,6)--(5,6)--(5,9)--(3,9)--(3,8)--(4,8)--(4,6);
\draw [line width = 1, fill=lightgray] (6,4)--(7,4)--(7,7)--(5,7)--(5,6)--(6,6)--(6,4);
\draw [line width = 1, fill=lightgray] (8,2)--(9,2)--(9,5)--(7,5)--(7,4)--(8,4)--(8,2);
\draw [line width = 1, fill=lightgray] (8,0)--(11,0)--(11,1)--(9,1)--(9,2)--(8,2)--(8,0);

\end{tikzpicture}
\caption{$11\times 11$ deficient square.}
\end{subfigure}
~~~~~~~~
\begin{subfigure}{0.3\textwidth}
\begin{tikzpicture}[scale = .3]
\draw (0,0) grid (13,13);

		\draw [line width = 1] (0,0)--(4,0)--(4,10)--(0,10)--(0,0);
		\draw [line width = 1] (0,2)--(4,2);
		\draw [line width = 1] (0,4)--(4,4);
		\draw [line width = 1] (0,6)--(4,6);
		\draw [line width = 1] (0,8)--(4,8);
		
		\draw [line width = 1] (4,0)--(8,0)--(8,6)--(4,6)--(4,0);
		\draw [line width = 1] (4,2)--(8,2);
		\draw [line width = 1] (4,4)--(8,4);
		
		\draw [line width = 1] (8,0)--(10,0)--(10,4)--(8,4)--(8,0);
		\draw [line width = 1] (11,1)--(13,1)--(13,5)--(11,5)--(11,1);
		
		\draw [line width = 1] (13,13)--(9,13)--(9,5)--(13,5)--(13,13);
		\draw [line width = 1] (13,11)--(9,11);
		\draw [line width = 1] (13,9)--(9,9);
		\draw [line width = 1] (13,7)--(9,7);
		\draw [line width = 1] (13,5)--(9,5);
		
		\draw [line width = 1] (5,7)--(9,7)--(9,13)--(5,13)--(5,7);
		\draw [line width = 1] (5,9)--(9,9);
		\draw [line width = 1] (5,11)--(9,11);
		
		\draw [line width = 1] (1,11)--(5,11)--(5,13)--(1,13)--(1,11);
		
		\draw [line width = 1, fill=lightgray] (13,0)--(13,1)--(11,1)--(11,2)--(10,2)--(10,0)--(13,0);
		\draw [line width = 1, fill=lightgray] (10, 2)--(11,2)--(11,5)--(9,5)--(9,4)--(10,4)--(10,2);
		\draw [line width = 1, fill=lightgray] (8, 4)--(9,4)--(9,7)--(7,7)--(7,6)--(8,6)--(8,4);
		\draw [line width = 1, fill=lightgray] (4,6)--(7,6)--(7,7)--(5,7)--(5,8)--(4,8)--(4,6);
		\draw [line width = 1, fill=lightgray] (4,9)--(5,9)--(5,11)--(2,11)--(2,10)--(4,10)--(4,9);
		\draw [line width = 1, fill=lightgray] (4,9)--(5,9)--(5,11)--(2,11)--(2,10)--(4,10)--(4,9);
		\draw [line width = 1, fill=lightgray] (0,10)--(2,10)--(2,11)--(1,11)--(1,13)--(0,13)--(0,10);    
		
		\draw [fill=darkgray] (4,8)--(5,8)--(5,9)--(4,9)--(4,8);
		
		\end{tikzpicture}
		\caption{$13\times 13$ deficient square.}
		\end{subfigure}
\caption{Tiling deficient squares.}
\label{fig:tiling_deficient_squares}
\end{figure}
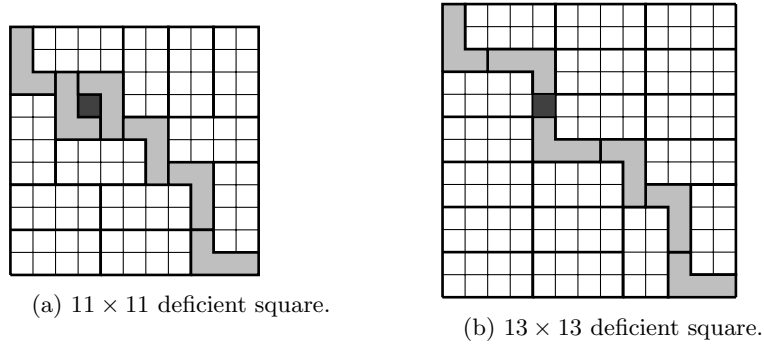

The tile sets $\mathcal{T}_4$ and $\mathcal{T}_4^+$ were introduced in \cite{CLNS} and \cite{N-replicating}, motivated by a problem in recreational mathematics. It is shown in \cite{CLNS} and \cite{N-replicating} that the tile sets have several remarkable  properties. Many of them are consequences of the fact that any tiling of the first quadrant by $\mathcal{T}_4$ or $\mathcal{T}_4^+$  follows the rectangular pattern, that is, the tiling reduces to a tiling by $2\times 4$ rectangles, in which every rectangle is tiled by two tiles from $\mathcal{T}_4$, and $2\times 2$ squares. This in turn, is a consequence of the fact that in any covering without overlaps of a region in the first quadrant bounded by a step 2 staircase and the coordinate axes, the $2$-squares are all covered by $2\times 4$ rectangles, covered by two tiles from $\mathcal{T}_4$, and $2\times 2$ squares with the coordinates of all vertices even. A bit unexpected, these results have applications to tilings of deficient rectangles. They also provide a natural mechanism for producing tilings with cracks.

Our main results are the following.

\begin{thm}\label{thm:main} Assume that a deficient rectangle is tiled by $\mathcal{T}_4$. Then the following are true:
\begin{enumerate}
\item The rectangle is a square of odd side $(2n+1)\times (2n+1), n\ge 2$. 
\item The missing cell has to be on the main NW--SE diagonal, in an odd position  if the square is $(4m+1)\times (4m+1)$ and in an even position if the square is $(4m+3)\times (4m+3)$. For all $m\ge 1$, all such positions give deficient squares that have tilings by 
$\mathcal{T}_4$.
\item The majority of the tiles (all but $n$ or $n+1$) in the tiling follow the rectangular pattern, that is, are paired and each pair tiles a $2\times 4$ rectangle. 
\item (Existence of the crack) The $n$ or $n+1$ tiles in an irregular position together with the missing cell form a NW--SE diagonal crack starting in the upper left corner of the square and ending in the lower right corner. 
\item (Location of the crack) The crack is located in a thin region symmetric about the diagonal, made out of a sequence of $3\times 3$ squares that overlap over one of the corner cells. 
\item (Symmetry of the crack) The crack divides the square in two parts of equal area. 
\item The number of possible cracks (counted twice if they allow for two different tilings) in a fixed deficient square is equal to $C_{2m}^m$ if the square is $(4m+1)\times (4m+1)$ and $2C_{2m}^m$ if the square is $(4m+3)\times (4m+3)$.
\item The number of tilings of a $(4m+1)\times (4m+1), m\ge 1,$ fixed deficient square by  $\mathcal{T}_4$ is equal to the number of tilings by dominoes of a $2m\times 2m$ square. 
\item The number of tilings of a $(4m+3)\times (4m+3), m\ge 1,$ fixed deficient square by  $\mathcal{T}_4$ is twice the number of tilings by dominoes of a $(2m+1)\times (2m+1)$ deficient square, with the missing cell placed on the main diagonal. 
\item The counting is realized by an explicit function that is a bijection in 8. and a double cover in 9., that takes a tiling by $\mathcal{T}_4$ into a tiling by dominoes. One takes a tiling of the deficient square by $\mathcal{T}_4$, eliminates the crack, replaces the pairs of tiles that cover $2\times 4$ rectangles by $2\times 4$ rectangles, reassembles the remaining two region in a square and then does a 1/2-homothety.
\item (Propagation of the crack) The crack and the tiling of a $(2n+1)\times (2n+1)$ deficient square can be extended (imbedded) into a crack (tiling) of a  $(2n+5)\times (2n+5)$ deficient square.
\end{enumerate}
\end{thm}

\begin{rmk} The tilings in Figure \ref{fig:tiling_deficient_squares} illustrate the statement of the theorem. The number of tilings in 8. is independent of the position of the missing cell on the diagonal and can be computed using Kasteleyn formula \cite{K}. The number of tilings in 9. depends on the position of the missing cell on the diagonal. For example, if $m=2$ the numbers of tilings by $\mathcal{T}_4$ of a $11\times 11$ deficient square are, in this order, 384, 224, 392, 224, 384. These are twice the numbers of tilings by dominoes of a $5\times 5$ deficient square.
\end{rmk}

\begin{thm}\label{thm:main2} Assume that a deficient rectangle is tiled by $\mathcal{T}_4^+$. Then the following are true:
\begin{enumerate}
\item The rectangle is a square of odd side $(2n+1)\times (2n+1), n\ge 2$. 
\item The missing cell has to be on the main NW--SE diagonal, in an odd position  if the square is $(4m+1)\times (4m+1)$ and in an even position if the square is $(4m+3)\times (4m+3)$. For all $m\ge 1$, all such positions give deficient squares that have tilings by $\mathcal{T}_4^+$.
\item The majority of the tiles (all but $n$ or $n+1$) in the tiling follow the rectangular pattern, that is, are either paired tetrominoes with each pair tiling a $2\times 4$ rectangle or $2\times 2$ squares. 
\item (Existence of the crack) The $n$ or $n+1$ tiles in an irregular position together with the missing cell form a NW--SE diagonal crack starting in the upper left corner of the square and ending in the lower right corner. 
\item (Location of the crack) The crack is located in a thin region symmetric about the diagonal, made out of a sequence of $3\times 3$ squares that overlap over one of the corner cells. 
\item The number of possible cracks (counted twice if they allow for two different tilings) in a fixed deficient square is $2^{2m}$ if the square is $(4m+1)\times (4m+1)$ and $2\cdot 2^{2m}$ if the square is $(4m+3)\times (4m+3)$.
\item The number of tilings of a $(4m+1)\times (4m+1), m\ge 1,$ deficient square by  $\mathcal{T}_4$ is equal to 
\begin{equation}
N(m)=\sum_{k=1}^{2m} 2^k N_k,
\end{equation}
where $N_k$ is the number of tilings by dominoes and monomers of a $2m\times 2m$ square that has the diagonal covered by exactly $k$ monomers (and $2m-k$ dominoes). 
\item The number of tilings of a $(4m+3)\times (4m+3), m\ge 1,$ deficient square by  $\mathcal{T}_4^+$ is $2N(m)$. 
\item The counting is realized by an explicit surjective function that takes a tiling by $\mathcal{T}_4$ into a tiling by dominoes and monomers. One takes a tiling of the deficient square by $\mathcal{T}_4$, eliminates the crack, replaces the pairs of tiles that cover $2\times 4$ rectangles by $2\times 4$ rectangles, reassembles the remaining two region in a square and then does a 1/2-homothety. If the image of a tiling by $\mathcal{T}_4^+$ has exactly $k$ monomers on the main diagonal, then the cardinality of the preimage of that image is $2^k$.
\item (Propagation of the crack) The crack and the tiling of a $(2n+1)\times (2n+1)$ square can be extended (imbedded) into a crack (tiling) of a  $(2n+5)\times (2n+5)$ square.
\end{enumerate}
\end{thm}

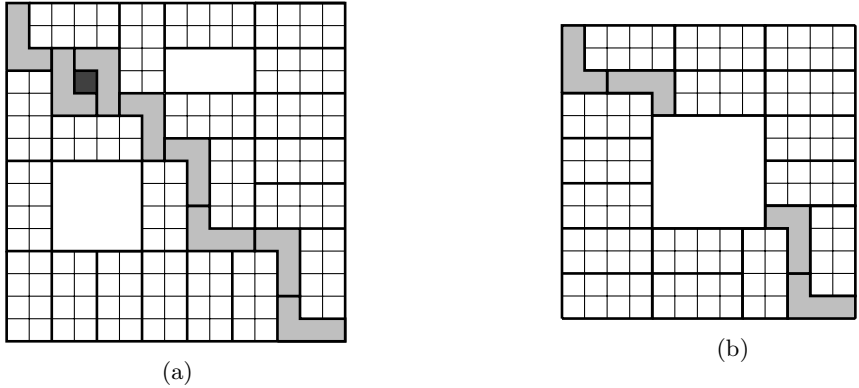
\begin{figure}[h]
\center
\begin{subfigure}{0.3\textwidth}
\begin{tikzpicture}[scale=.3]
\draw (0,-4) grid (15,11);
\draw [fill = darkgray] (3,7)--(4,7)--(4,8)--(3,8)--(3,7);

\draw [line width = 1] (0,0)--(0,-4)--(15,-4)--(15,11)--(11,11);
\draw [line width = 1] (4,0)--(0, 0) -- (0, 4);
\draw [line width = 1] (4,0)-- (8,0);

\draw [line width = 1] (11,1)--(11,5);
\draw [line width = 1] (11,5)--(11,11)--(5,11);
\draw [line width = 1] (1,11)--(5,11);
\draw [line width = 1] (0,8)--(0,4);

\draw [line width = 1, fill=lightgray] (0,8)--(2,8)--(2,9)--(1,9)--(1,11)--(0,11)--(0,8);
\draw [line width = 1, fill=lightgray] (2,6)--(4,6)--(4,7)--(3,7)--(3,9)--(2,9)--(2,6);
\draw [line width = 1, fill=lightgray] (4,6)--(5,6)--(5,9)--(3,9)--(3,8)--(4,8)--(4,6);
\draw [line width = 1, fill=lightgray] (6,4)--(7,4)--(7,7)--(5,7)--(5,6)--(6,6)--(6,4);
\draw [line width = 1, fill=lightgray] (8,2)--(9,2)--(9,5)--(7,5)--(7,4)--(8,4)--(8,2);
\draw [line width = 1, fill=lightgray] (8,0)--(11,0)--(11,1)--(9,1)--(9,2)--(8,2)--(8,0);

\draw [line width = 1, fill=lightgray] (12,-2)--(13,-2)--(13,1)--(11,1)--(11,0)--(12,0)--(12,-2);
\draw [line width = 1, fill=lightgray] (12,-4)--(15,-4)--(15,-3)--(13,-3)--(13,-2)--(12,-2)--(12,-4);

\draw [line width = 1] (2,0)--(2,-4);
\draw [line width = 1] (4,0)--(4,-4);
\draw [line width = 1] (6,0)--(6,-4);
\draw [line width = 1] (8,0)--(8,-4);
\draw [line width = 1] (10,0)--(10,-4);

\draw [line width = 1] (15,1)--(13,1);
\draw [line width = 1] (15,3)--(11,3);
\draw [line width = 1] (15,5)--(11,5);
\draw [line width = 1] (15,7)--(11,7);
\draw [line width = 1] (15,9)--(11,9);
\draw [line width = 1] (15,11)--(11,11);
\draw [line width = 1] (0,4)--(2,4)--(2,6);
\draw [line width = 1] (5,9)--(5,11);
\draw [line width = 1] (7,9)--(7,11);
\draw [line width = 1] (9,5)--(11,5);

\draw [line width=1, fill=white] (2,0)--(6,0)--(6,4)--(2,4)--(2,0);
\draw [line width=1, fill=white] (7,7)--(11,7)--(11,9)--(7,9)--(7,7);

\end{tikzpicture}
\caption{}
\end{subfigure}
~~~~~~~~~~~~~~~~~~~~~~
\begin{subfigure}{0.3\textwidth}
\begin{tikzpicture}[scale = .3]
\draw (0,0) grid (13,13);

		\draw [line width = 1] (0,0)--(4,0)--(4,10)--(0,10)--(0,0);
		\draw [line width = 1] (0,2)--(4,2);
		\draw [line width = 1] (0,4)--(4,4);
		\draw [line width = 1] (0,6)--(4,6);
		\draw [line width = 1] (0,8)--(4,8);
		
		\draw [line width = 1] (4,0)--(8,0)--(8,6)--(4,6)--(4,0);
		\draw [line width = 1] (4,2)--(8,2);
		\draw [line width = 1] (4,4)--(8,4);
		
		\draw [line width = 1] (8,0)--(10,0)--(10,4)--(8,4)--(8,0);
		\draw [line width = 1] (11,1)--(13,1)--(13,5)--(11,5)--(11,1);
		
		\draw [line width = 1] (13,13)--(9,13)--(9,5)--(13,5)--(13,13);
		\draw [line width = 1] (13,11)--(9,11);
		\draw [line width = 1] (13,9)--(9,9);
		\draw [line width = 1] (13,7)--(9,7);
		\draw [line width = 1] (13,5)--(9,5);
		
		\draw [line width = 1] (5,7)--(9,7)--(9,13)--(5,13)--(5,7);
		\draw [line width = 1] (5,9)--(9,9);
		\draw [line width = 1] (5,11)--(9,11);
		
		\draw [line width = 1] (1,11)--(5,11)--(5,13)--(1,13)--(1,11);
		
		\draw [line width = 1, fill=lightgray] (13,0)--(13,1)--(11,1)--(11,2)--(10,2)--(10,0)--(13,0);
		\draw [line width = 1, fill=lightgray] (10, 2)--(11,2)--(11,5)--(9,5)--(9,4)--(10,4)--(10,2);
		\draw [line width = 1, fill=lightgray] (8, 4)--(9,4)--(9,7)--(7,7)--(7,6)--(8,6)--(8,4);
		\draw [line width = 1, fill=lightgray] (4,6)--(7,6)--(7,7)--(5,7)--(5,8)--(4,8)--(4,6);
		\draw [line width = 1, fill=lightgray] (4,9)--(5,9)--(5,11)--(2,11)--(2,10)--(4,10)--(4,9);
		\draw [line width = 1, fill=lightgray] (4,9)--(5,9)--(5,11)--(2,11)--(2,10)--(4,10)--(4,9);
		\draw [line width = 1, fill=lightgray] (0,10)--(2,10)--(2,11)--(1,11)--(1,13)--(0,13)--(0,10);    
		
		\draw [fill=darkgray] (4,8)--(5,8)--(5,9)--(4,9)--(4,8);
		
		\draw [line width = 1, fill=white] (9,4)--(9,9)--(4,9)--(4,4)--(9,4);
		
		\end{tikzpicture}
		\caption{}
		\end{subfigure}
\caption{Tiling irregular regions.}
\label{fig:tiling_deficient_regions}
\end{figure}

Our results about the existence and the properties of tilings with cracks can be easily extended to more irregular regions, and even to regions with extra holes. The idea of the proof is to complete these regions to a full deficient square and then apply the results in Theorems \ref{thm:main} and \ref{thm:main2}. Some typical examples are shown in Figure \ref{fig:tiling_deficient_regions}. 

The results in this note can be generalized to tilings of deficient rectangles by the sets of ribbon $L$ $n$-ominoes, $n$ even. These tile sets are introduced and studied in \cite{nitica-L-shaped}. We will do this in a future paper. Probably similar results may also be derived for more general tile sets appearing from dissection of rectangles that follow the rectangular pattern. These tile sets are described in \cite{CFNS} and \cite{N-last}.

Tilings of deficient rectangles by $L$-tetrominoes, with all 8 symmetries allowed, are described in \cite{N-L-shaped}.  

\section{Proof of Theorem \ref{thm:main}}

\begin{figure}[h]
\center
\begin{tikzpicture}[scale=.3]
\draw (0,0) grid (17,11);
\draw [line width = 1] (0,0)--(17,0)--(17,11)--(0,11)--(0,0);

\draw [line width = 1] (2,0)--(2,10);
\draw [line width = 1] (4,0)--(4, 8);
\draw [line width = 1] (6,0)-- (6,6);
\draw [line width = 1] (8,0)--(8,4);
\draw [line width = 1] (0,8)--(0,4);

\draw [line width = 1] (0,10)--(2,10);
\draw [line width = 1] (0,8)--(4,8);
\draw [line width = 1] (0,6)--(6,6);
\draw [line width = 1] (0,4)--(8,4);
\draw [line width = 1] (0,2)--(10,2)--(10,0);

\draw [line width = 1] (15,11)--(15,1);
\draw [line width = 1] (13,11)--(13,3);
\draw [line width = 1] (11,11)--(11,5);
\draw [line width = 1] (9,11)--(9,7);
\draw [line width = 1] (7,11)--(7,9);
\draw [line width = 1] (7,9)--(17,9);
\draw [line width = 1] (9,7)--(17,7);
\draw [line width = 1] (11,5)--(17,5);
\draw [line width = 1] (13,3)--(17,3);
\draw [line width = 1] (15,1)--(17,1);

\draw [line width = 1, fill=lightgray] (0,10)--(2,10)--(2,9)--(3,9)--(3,11)--(0,11)--(0,10);
\draw [line width = 1, fill=lightgray] (2,8)--(4,8)--(4,7)--(5,7)--(5,9)--(2,9)--(2,8);
\draw [line width = 1, fill=lightgray] (4,6)--(6,6)--(6,5)--(7,5)--(7,7)--(4,7)--(4,6);
\draw [line width = 1, fill=lightgray] (6,4)--(8,4)--(8,3)--(9,3)--(9,5)--(6,5)--(6,4);
\draw [line width = 1, fill=lightgray] (8,2)--(10,2)--(10,1)--(11,1)--(11,3)--(8,3)--(8,2);

\node at (10.5,.5) {*};
\end{tikzpicture}
\caption{A general rectangle.}
\label{fig:gen-rectangle}
\end{figure}
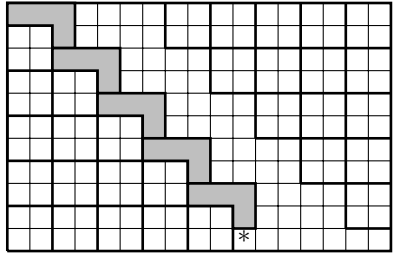

It follows from \cite{CLNS} that any tiling by $\mathcal{T}_4$ of a region in the first quadrant bounded by the coordinate axes and a staircase of even origin and step 2 has to follow the rectangular pattern. For simplicity we will call such a region \emph{staircase}. Inside any $(4p+1)\times (4q+1)$ or  $(4p+3)\times (4q+3)$ rectangle we can fit two maximal regions as above. See Figure \ref{fig:gen-rectangle}. Assume that the missing cell is not placed on the diagonal adjacent to one of the staircase regions or inside that region. This forces the appearance of the gray staircase which ends in the cell * that cannot be tiled. We conclude that the deficient rectangle has to be a square and the missing cell has to be placed on the diagonal.

\begin{figure}[h]
\center
\begin{tikzpicture}[scale=.3]
\draw (0,0) grid (11,11);
\draw [line width = 1] (0,0)--(11,0)--(11,11)--(0,11)--(0,0);

\draw [line width = 1] (2,0)--(2,8);
\draw [line width = 1] (4,0)--(4, 6);
\draw [line width = 1] (6,0)-- (6,4);
\draw [line width = 1] (8,0)--(8,2);

\draw [line width = 1] (0,8)--(2,8);
\draw [line width = 1] (0,6)--(4,6);
\draw [line width = 1] (0,4)--(6,4);
\draw [line width = 1] (0,2)--(8,2)--(8,0);

\draw [line width = 1] (9,11)--(9,3);
\draw [line width = 1] (7,11)--(7,5);
\draw [line width = 1] (5,11)--(5,7);
\draw [line width = 1] (3,11)--(3,9);

\draw [line width = 1] (3,9)--(11,9);
\draw [line width = 1] (5,7)--(11,7);
\draw [line width = 1] (7,5)--(11,5);
\draw [line width = 1] (9,3)--(11,3);
\end{tikzpicture}
\caption{The central region.}
\label{fig:central_region}
\end{figure}
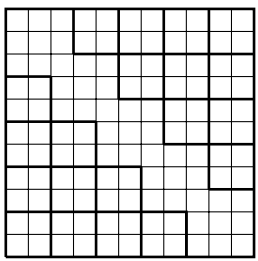

Assume now that we have a deficient square. See Figure \ref{fig:central_region}. It follows from \cite{CLNS} that all $2\times 2$ marked squares are tiled following the rectangular pattern. We study now the tiling of the remaining central region consisting of a sequence of $3\times 3$ squares centered about the main diagonal. We observe first that some of the central region has to be covered by $2\times 2$ squares that are parts of $2\times 4$ rectangles originating in the region covered by the staircases. Indeed, doing a chessboard coloring (say black, white) of the $2$-squares in the staircase region, we need to have the same number of black and white squares. This is due to the fact that any $2\times 4$ rectangle is covered by a black and a white $2$-square. If the deficient square is $(4m+1)\times (4m+1)$ then we have a deficiency of $m$ $2$-squares for each of the staircases and if the deficient square is $(4m+3)\times (4m+3)$ then we have a deficiency of $m$ $2$-squares for each of the staircases. In the first case this forces all $3\times 3$ squares in the central region to contain a $2$-square and in the second case forces all but one of the $3\times 3$ squares in the central region to contain a $2$-square. The $2$-squares are distributed evenly between the two staircases. The way in which we partition them in two equal parts, as we see below, can be arbitrary. 

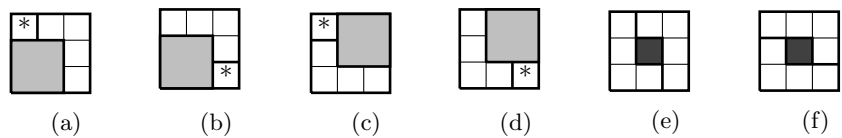
\begin{figure}[h]
\center
\begin{subfigure}{0.1\textwidth}
\begin{tikzpicture}[scale=.35]
\draw (0,0) grid (3,3);

\draw [line width = 1] (0,0)--(3,0)--(3,3)--(0,3)--(0,0);

\draw [line width = 1, fill=lightgray] (0,0)--(2, 0)--(2, 2)--(0,2)--(0,0);

\draw [line width = 1] (1,2)--(1,3);
\node at (.5,2.5) {*};

\end{tikzpicture}
\caption{}
\end{subfigure}
~~
\begin{subfigure}{0.1\textwidth}
\begin{tikzpicture}[scale=.35]
\draw (0,0) grid (3,3);

\draw [line width = 1] (0,0)--(3,0)--(3,3)--(0,3)--(0,0);

\draw [line width = 1, fill=lightgray] (0,0)--(2, 0)--(2, 2)--(0,2)--(0,0);

\draw [line width = 1] (2,1)--(3,1);
\node at (2.5,.5) {*};

\end{tikzpicture}
\caption{}
\end{subfigure}
~~
\begin{subfigure}{0.1\textwidth}
\begin{tikzpicture}[scale=.35]
\draw (0,0) grid (3,3);

\draw [line width = 1] (0,0)--(3,0)--(3,3)--(0,3)--(0,0);

\draw [line width = 1, fill=lightgray] (1,1)--(3, 1)--(3, 3)--(1,3)--(1,1);

\draw [line width = 1] (0,2)--(1,2);
\node at (.5,2.5) {*};

\end{tikzpicture}
\caption{}
\end{subfigure}
~~
\begin{subfigure}{0.1\textwidth}
\begin{tikzpicture}[scale=.35]
\draw (0,0) grid (3,3);

\draw [line width = 1] (0,0)--(3,0)--(3,3)--(0,3)--(0,0);

\draw [line width = 1, fill=lightgray] (1,1)--(3, 1)--(3, 3)--(1,3)--(1,1);

\draw [line width = 1] (2,0)--(2,1);
\node at (2.5,.5) {*};

\end{tikzpicture}
\caption{}
\end{subfigure}
~~
\begin{subfigure}{0.1\textwidth}
\begin{tikzpicture}[scale=.35]
\draw (0,0) grid (3,3);

\draw [line width = 1] (0,0)--(3,0)--(3,3)--(0,3)--(0,0);

\draw [line width = 1, fill=darkgray] (1,1)--(2, 1)--(2, 2)--(1,2)--(1,1);

\draw [line width = 1] (1,2)--(1,3);
\draw [line width = 1] (2,0)--(2,1);

\end{tikzpicture}
\caption{}
\end{subfigure}
~~
\begin{subfigure}{0.1\textwidth}
\begin{tikzpicture}[scale=.35]
\draw (0,0) grid (3,3);

\draw [line width = 1] (0,0)--(3,0)--(3,3)--(0,3)--(0,0);

\draw [line width = 1, fill=darkgray] (1,1)--(2, 1)--(2, 2)--(1,2)--(1,1);

\draw [line width = 1] (2,1)--(3,1);
\draw [line width = 1] (0,2)--(1,2);

\end{tikzpicture}
\caption{}
\end{subfigure}
\caption{Tilings of a $3\times 3$ square.}
\label{fig:three}
\end{figure}

We study now the region left uncovered in the central region after the addition of the $2$-squares. The possible tilings of a $3\times 3$ square are shown in Figure \ref{fig:three}. In cases a) through d) the cell labeled * has to be covered by the missing cell in the deficient square or by a cell that is part of a tile in $\mathcal{T}_4$ originating in a different $3\times 3$ square. In  cases e) and f) the central cell in $3\times 3$ square has to be the missing cell. If the deficient square is of type $(4m+1)\times (4m+1)$ then the $3\times 3$ squares are covered by the cases a) through d). The missing cell has to be in an odd position on the main diagonal covering the corner of a $3\times 3$ square. If the deficient square is of type $(4m+3)\times (4m+3)$ then all but one of the $3\times 3$ squares are covered by the cases a) through d) and one of them is covered by the case e) or f). The missing cell has to be in an even position on the main diagonal covering the center of a $3\times 3$ square. It is easy to see that once the $2$-squares and the missing cell are places in the central region, the rest of the crack can be tiled by $\mathcal{T}_4$ in a unique way for $(4m+1)\times (4m+1)$ squares and in two ways for $(4m+3)\times (4m+3)$ squares.

\begin{figure}[h]
\center
\begin{tikzpicture}[scale=.3]
\draw (0,0) grid (17,17);
\draw [line width = 1,<->] (0,17.5)--(0,0)--(17.5,0);

\draw [line width = 1] (2,0)--(2,14);
\draw [line width = 1] (4,0)--(4, 12);
\draw [line width = 1] (6,0)-- (6,10);
\draw [line width = 1] (8,0)--(8,8);
\draw [line width = 1] (10,0)--(10,6);
\draw [line width = 1] (12,0)--(12,4);
\draw [line width = 1, fill=lightgray] (0,14)--(2,14)--(2,16)--(0,16)--(0,14);
\draw [line width = 1, fill=lightgray] (2,12)--(4,12)--(4,14)--(2,14)--(2,12);
\draw [line width = 1, fill=lightgray] (4,10)--(6,10)--(6,12)--(4,12)--(4,10);
\draw [line width = 1, fill=lightgray] (14,0)--(14,2)--(16,2)--(16,0)--(14,0);

\draw [line width = 1] (0,14)--(2,14);
\draw [line width = 1] (0,12)--(4,12);
\draw [line width = 1] (0,10)--(6,10);
\draw [line width = 1] (0,8)--(8,8);
\draw [line width = 1] (0,6)--(10,6);
\draw [line width = 1] (0,4)--(12,4);
\draw [line width = 1] (0,2)--(14,2)--(14,0);

\draw [line width = 1] (20,2)--(32,2)--(32,0);
\draw [line width = 1] (20,4)--(30,4);
\draw [line width = 1] (20,6)--(28,6);
\draw [line width = 1] (20,8)--(26,8);
\draw [line width = 1] (20,10)--(24,10);
\draw [line width = 1] (20,12)--(22,12);

\draw [line width = 1] (22,0)--(22,12);
\draw [line width = 1] (24,0)--(24,10);
\draw [line width = 1] (26,0)--(26,8);
\draw [line width = 1] (28,0)--(28,6);
\draw [line width = 1] (30,0)--(30,4);
\draw [line width = 1] (32,0)--(32,2);

\draw [line width = 4] (13,1)--(15,1);
\draw [line width = 4] (1,13)--(1,15);
\draw [line width = 4] (3,11)--(3,13);
\draw [line width = 4] (5,9)--(5,11);

\draw [line width = 5, ->] (17.5,10)--(19.5,10);

\draw (20,0) grid (35,17);
\draw [line width = 1,<->] (20,17.5)--(20,0)--(35.5,0);

\draw [line width = 1, fill=lightgray] (30,2)--(32,2)--(32,4)--(30,4)--(30,2);
\draw [line width = 1, fill=lightgray] (28,4)--(30,4)--(30,6)--(28,6)--(28,4);
\draw [line width = 1, fill=lightgray] (26,6)--(28,6)--(28,8)--(26,8)--(26,6);

\end{tikzpicture}
\caption{The induction step.}
\label{fig:addition}
\end{figure}
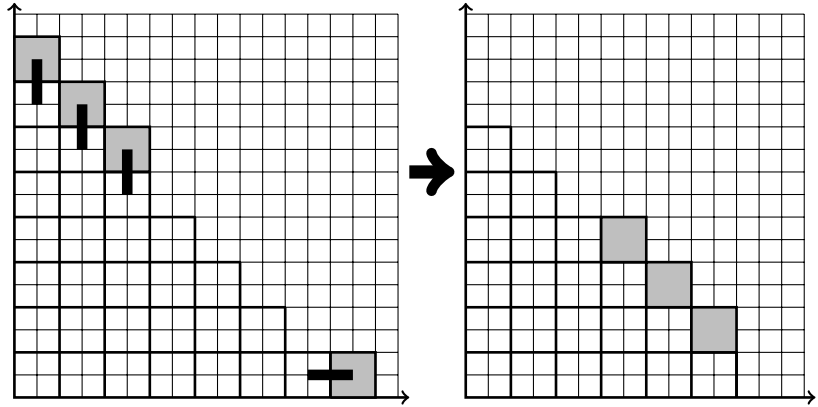

We show now that after the addition of the $2$-squares to the staircase the resulting region can be tiled by $2\times 4$ rectangles. This is obvious if the number of rows in the staircase is 1 or 2. If the number of rows is larger, we do induction on the number of rows. The induction step is illustrated in Figure \ref{fig:addition} and consists in removing $2\times 4$ rectangles containing the additional $2$-squares. We start removing $2\times 4$ rectangles from the top of the staircase. Each removal deletes also a $2$-square on the largest diagonal of the staircase. In order to avoid ambiguity, we want the $2$-squares on the diagonal to be eliminated in order, from the top to the bottom. What is left is a smaller staircase with the additional $2$-squares added which can be tiled by $2\times 4$ rectangles due to the induction hypothesis.

To finish the proof of the counting results we need to show that any tiling by dominoes of a $2m\times 2m$ square or of a $(2m+1)\times (2m+1)$ deficient square with the missing cell on the main diagonal  can be divided in tilings of complementary staircases regions with extra squares added. Observe that any of the $2m$ available cells on the diagonal is covered by a different domino. To obtain the splitting, divide the set of dominoes that cover the cells on the diagonal in two equal parts and assign them to the lower, respectively upper, maximal staircase in the original $2m\times 2m$ or $(2m+1)\times (2m+1)$ square. This argument also allows to count the number of cracks. The crack is uniquely determined by the partition of extra $2\times 2$ squares in two equal parts between the upper and lower staircases. Independent of the size of the deficient square, this partition can be done in $C_{2m}^m$ ways.

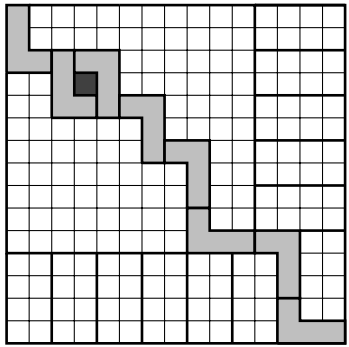
\begin{figure}[h]
\center
\begin{tikzpicture}[scale=.3]
\draw (0,-4) grid (15,11);
\draw [fill = darkgray] (3,7)--(4,7)--(4,8)--(3,8)--(3,7);

\draw [line width = 1] (0,0)--(0,-4)--(15,-4)--(15,11)--(11,11);
\draw [line width = 1] (4,0)--(0, 0) -- (0, 4);
\draw [line width = 1] (4,0)-- (8,0);
\draw [line width = 1] (11,1)--(11,5);
\draw [line width = 1] (11,5)--(11,11)--(5,11);
\draw [line width = 1] (1,11)--(5,11);
\draw [line width = 1] (0,8)--(0,4);

\draw [line width = 1, fill=lightgray] (0,8)--(2,8)--(2,9)--(1,9)--(1,11)--(0,11)--(0,8);
\draw [line width = 1, fill=lightgray] (2,6)--(4,6)--(4,7)--(3,7)--(3,9)--(2,9)--(2,6);
\draw [line width = 1, fill=lightgray] (4,6)--(5,6)--(5,9)--(3,9)--(3,8)--(4,8)--(4,6);
\draw [line width = 1, fill=lightgray] (6,4)--(7,4)--(7,7)--(5,7)--(5,6)--(6,6)--(6,4);
\draw [line width = 1, fill=lightgray] (8,2)--(9,2)--(9,5)--(7,5)--(7,4)--(8,4)--(8,2);
\draw [line width = 1, fill=lightgray] (8,0)--(11,0)--(11,1)--(9,1)--(9,2)--(8,2)--(8,0);

\draw [line width = 1, fill=lightgray] (12,-2)--(13,-2)--(13,1)--(11,1)--(11,0)--(12,0)--(12,-2);
\draw [line width = 1, fill=lightgray] (12,-4)--(15,-4)--(15,-3)--(13,-3)--(13,-2)--(12,-2)--(12,-4);

\draw [line width = 1] (2,0)--(2,-4);
\draw [line width = 1] (4,0)--(4,-4);
\draw [line width = 1] (6,0)--(6,-4);
\draw [line width = 1] (8,0)--(8,-4);
\draw [line width = 1] (10,0)--(10,-4);

\draw [line width = 1] (15,1)--(13,1);
\draw [line width = 1] (15,3)--(11,3);
\draw [line width = 1] (15,5)--(11,5);
\draw [line width = 1] (15,7)--(11,7);
\draw [line width = 1] (15,9)--(11,9);
\draw [line width = 1] (15,11)--(11,11);

\end{tikzpicture}
\caption{Propagation of the crack.}
\label{fig:propagation}
\end{figure}

The propagation of the crack is illustrated in Figure \ref{fig:propagation} and it is self explanatory. 

$\square$

\section{Proof of Theorem \ref{thm:main2}}

The proof of Theorem \ref{thm:main2} is similar to that of Theorem \ref{thm:main}. Figure \ref{fig:three} shows that the appearance of the extra $2\times 2$ square inside the $3\times 3$ squares is forced and the only possible coverings of a $3\times 3$ square are shown in that figure. The differences in counting appear due to the fact that the extra $2\times 2$ squares that are added in the central region do not need to be divided in two equal parts, but rather can be divided arbitrarily. If the extra $2\times 2$ squares are placed, then the region covered by the crack is well defined. As a $2\times 2$ square can be placed in only two positions, the number of cracks is $2^{2m}$ for a $(4m+1)\times (4m+1)$ square and $2\cdot 2^{2m}$ for a $(4m+3)\times (4m+3)$ square. The extra factor of $2$ appears due to the fact that in the last case the region supporting the crack can be tiled in two ways. See Figure \ref{fig:three}, e), f). Same arguments justify the extra factors of $2^k$ and $2$ appearing in the statements 7. and 8. in Theorem \ref{thm:main2}.

\section*{Bibliography}
		
\bibliographystyle{plain}

\begin{thebibliography}{9999}

\bibitem{CFNS} A.~Calderon, S.~Fairchild, V.~Nitica, S.~Simon, Tilings of quadrants by $L$-ominoes and notched rectangles, Topics in Recreational Mathematics, {\em 7}, 2016, pp. 39--75.

\bibitem{CLNS} M.~Chao, D.~Levenstein, V.~Nitica, R.~Sharp, A coloring invariant for ribbon $L$-tetrominoes, Discrete Mathematics, {\em 313} (2013) 611--621.

\bibitem{K} P.~M.~Kasteleyn, The statistics of dimers on a lattice, Physica {\em 27}, (1961) 1209--1225.

\bibitem{N-L-shaped} C.~Nitica, V.~Nitica, Tiling a deficient rectangle by $L$-tetrominoes, Journal of Recreational Mathematics  {\em 33} (4) (2004-2005) 259--271.

\bibitem{N-replicating} V.~Nitica, Rep-tiles revisited, in the volume {\it MASS Selecta: Teaching and Learning
Advanced Undergraduate Mathematics}, American Mathematical Society, 2003

\bibitem{nitica-L-shaped} V.~Nitica,  Every tiling of the first quadrant by ribbon $L$ $n$-ominoes follows the rectangular
pattern. Open Journal of Discrete Mathematics, {\em 5}, (2015) 11--25.

\bibitem{N-last} V.~Nitica, On tilings of quadrants and rectangles and rectangular pattern, Open Journal of Discrete Mathematics {\em 6}, (2016) 351--371.

\end{thebibliography}

\end{document}